\documentclass{elsart}
\usepackage{amsmath,amssymb}
\usepackage[english,francais]{babel}

\newcommand{\loc}{\mathrm{loc}}
\renewcommand{\div}{\mathrm{div}}
\newcommand{\er}{\mathbf{e}_\rho}

\newtheorem{theorem}{Theorem}[section]

\newtheorem{e-proposition}[theorem]{Proposition}
\newtheorem{corollary}[theorem]{Corollary}
\newtheorem{e-definition}[theorem]{Definition\rm}
\newtheorem{remark}{\it Remark\/}

\newcommand{\QED}{\hfill\ensuremath{\square}}

\def\og{\leavevmode\raise.3ex\hbox{$\scriptscriptstyle\langle\!\langle$~}}
\def\fg{\leavevmode\raise.3ex\hbox{~$\!\scriptscriptstyle\,\rangle\!\rangle$}}

\journal{}

\begin{document}

\centerline{}

\begin{frontmatter}

\selectlanguage{english}
\title{A generalization of the Hopf-Cole transformation for stationary Mean Field Games systems}
\selectlanguage{english}
\author{Marco Cirant},
\ead{marco.cirant@unimi.it}

\address{Dipartimento di Matematica ``F. Enriques''\\ Universit\`a di Milano\\ Via C. Saldini, 50\\ 20133--Milano, Italy}

\thanks{The author is supported by a Post-Doc Fellowship from the Universit\`a degli Studi di Milano.}

\begin{abstract}
In this note we propose a transformation which decouples stationary Mean Field Games systems with superlinear Hamiltonians of the form $|p|^{r'}$, $r' > 1$, and turns the Hamilton-Jacobi-Bellman equation into a quasi-linear equation involving the $r$-Laplace operator. Such a transformation requires an assumption on solutions of the system, which is satisfied for example in space dimension one or if solutions are radial.

\vskip 0.5\baselineskip


\end{abstract}


\end{frontmatter}

\section{Introduction}

Mean Field Games (briefly MFG) is a branch of Dynamic Games which has been proposed independently by Lasry, Lions \cite{jeux1}, \cite{jeux2}, \cite{LasryLions}, \cite{LionsVideo} and Caines, Huang, Malham\'e \cite{HCM}, and aims at modeling and analyzing decision processes involving a very large number of indistinguishable rational agents. In MFG, every agent belonging to a population of infinite individuals has the goal of minimizing some cost which depends on his own state and on the average distribution of the other players. Suppose that the state of a typical player is driven by the stochastic differential equation
\[
d X_s = -\alpha_s ds + \sqrt{2 \nu} \, d B_s \in \Omega,
\]
where $\alpha_s$ is the control, $B_s$ is a Brownian motion, $\nu > 0$ and the domain $\Omega \subseteq \Rset^d$, $d \ge 1$, is the so-called state space. Suppose also that the cost functional has the long-time-average form
\[
\mathcal{J}(X_0, \alpha) = \liminf_{T \rightarrow \infty} \frac{1}{T} \int_0^T \mathbb{E} [L(\alpha_s) + f(X_s, \hat{m}_s)] ds,
\]
where the (convex) Lagrangian function $L(\alpha)$ is associated to the cost paid by the player to change his own state, and the term involving $f$, that we assume to be a $C^1(\Omega \times [0, \infty))$ function, is the cost paid for being at state $x \in \Omega$, and it depends on the empirical density $\hat{m}_s$ of the other players. Then, under the assumption that players are indistinguishable, it has been shown that equilibria of the game (in the sense of Nash) are captured by the following system of non-linear elliptic equations
\begin{equation}\label{MFG}
\left\{
\begin{array}{rll}
\text{(HJB)} \quad & - \nu \Delta u(x) + H(Du(x)) + \lambda = f(x,m(x)) & \quad \text{in $\Omega$} \\
\text{(K)} \quad & - \nu \Delta m(x) - \div(D H(Du(x)) \, m(x)) = 0  &\quad \text{in $\Omega$}  \\
& \int_\Omega m(x) dx = 1, \, m \ge 0.
\end{array}
\right.
\end{equation}
Here, $H$ denotes the Legendre transform of $L$. The two unknowns $u, \lambda$ in the Hamilton-Jacobi-Bellman equation \eqref{MFG}-(HJB) provide respectively the optimal control of a typical player, given in feedback form by $\alpha^* : x \mapsto -D H(Du(x))$, and the average cost $\mathcal{J}(X_0, \alpha^*)$. On the other hand, the solution $m$ of the Kolmogorov equation \eqref{MFG}-(K) is the stationary distribution of players implementing the optimal strategy, that is the long time behavior of the whole population playing in an optimal way. The two equations are coupled via the cost function $f$. Note that a set of boundary conditions is usually associated to \eqref{MFG}, for example $u,m$ can be assumed to be periodic if $\Omega = (0,1)^d$ or Neumann conditions are imposed if $\Omega$ is bounded and $X_s$ is subject to reflection at $\partial \Omega$. In some models $\Omega$ is the whole $\Rset^d$.

A relevant class of MFG models assumes that the Lagrangian function $L$ has the form $L(\alpha) = \frac{l_0}{r} |\alpha|^{r}$, $l_0 > 0 , r > 1$. Consequently, the Hamiltonian function $H$ becomes
\begin{equation}\label{H}
H(p) = \frac{h_0}{r'} |p|^{r'}, \quad h_0 = l_0^{1-r'} > 0 , r' = \frac{r}{r-1} > 1.
\end{equation}
In the particular situation where $L$ and $H$ are \textit{quadratic} (namely $r= r' = 2$), it has been pointed out (see \cite{jeux1}, \cite{LionsVideo}) that the so-called Hopf-Cole transformation decouples \eqref{MFG}, and reduces it to a single elliptic semilinear equation of generalized Hartree type. Precisely, let $\varphi := c e^{-\frac{u}{2\nu}}$, $c>0$, then \eqref{MFG}-(HJB) reads (setting for simplicity $h_0 = 1$)
\begin{equation}\label{semilin}
- 2\nu^2 \Delta \varphi + (f(x, m) - \lambda) \varphi = 0 \quad \text{in $\Omega$} \\
\end{equation}
for all $c > 0$. Moreover, if we set $c^2 = \left(\int_\Omega e^{-\frac{u}{\nu}} \right)^{-1}$, an easy computation shows that $\varphi^2$ is also a solution of \eqref{MFG}-(K), so if uniqueness of solutions for such equation holds (that is true, for example, if suitable boundary conditions are imposed), then $m = \varphi^2$, and therefore $\varphi$ becomes the only unknown in \eqref{semilin}.
This transformation can be exploited to study quadratic MFG systems, both from the theoretical and numerical point of view. This strategy is adopted, for example, in the works \cite{cardalialong}, \cite{gomes}, \cite{gueant09}, \cite{gueant12}.

The aim of this note is to show that if $r' \neq 2$ there exists a similar change of variables, involving a suitable power of $m$, that in some cases decouples \eqref{MFG} and turns the Hamilton-Jacobi-Bellman equation \eqref{MFG}-(HJB) into a quasi-linear equation of the form
\begin{equation}\label{rlapl}
\left\{
\begin{array}{ll}
-\mu \Delta_r \varphi + (f(x, \varphi^r)- \lambda) \varphi^{r-1} = 0 & \text{in $\Omega$}, \\
\int_\Omega \varphi^r dx = 1, \, \varphi > 0, \, \mu = \nu \left(\frac{\nu r}{h_0}\right)^{r-1},
\end{array}
\right.
\end{equation}
where $\Delta_r \varphi = \div(|D\varphi|^{r-2} D\varphi)$ is the standard $r$-Laplace operator ($r$ is the conjugate exponent of $r'$). Such a transformation is in particular possible if the vector field $\nu D m + D H(Du) m$, which is divergence-free because of \eqref{MFG}-(K), is identically zero on $\Omega$. 

Let us briefly recall in which sense $(u,m,\lambda)$ solves \eqref{MFG}.
\begin{e-definition} We say that a triple $(u,m,\lambda) \in C^2({\Omega}) \times W^{1,2}_\loc(\Omega) \times \Rset$ is a (local) solution of \eqref{MFG} if $u,\lambda$ solve pointwise \eqref{MFG}-\textnormal{(HJB)} and $m$ solves \eqref{MFG}-\textnormal{(K)} in the weak sense, namely
\begin{equation}\label{weakolmo}
\nu \int_\Omega Dm \cdot D\xi + \int_\Omega m \, D H(Du) \cdot D \xi = 0 \quad \forall \xi \in C^{\infty}_0(\Omega).
\end{equation}
We say that a couple $(\varphi,\lambda) \in (W^{1,r}_\loc({\Omega}) \cap L^{\infty}_\loc({\Omega}))\times \Rset$ is a solution of \eqref{rlapl} if
\[
\mu \int_\Omega |D \varphi|^{r-2} D \varphi \cdot D \xi +  \int_\Omega (f(x, \varphi^r)- \lambda) \varphi^{r-1} \xi = 0 \quad \forall \xi \in C^{\infty}_0(\Omega).
\]
\end{e-definition}

Then, the transformation can be stated as follows.
\begin{theorem}\label{genHCp} Suppose that $H$ satisfies \eqref{H}.
\begin{itemize}
\item[\textbf{a)}] Let $(u,m,\lambda)$ be a solution of \eqref{MFG}. If the following equality holds,
\begin{equation}\label{crucialeq}
\nu Dm + h_0 m |Du|^{r'-2}Du = 0 \quad \text{a.e. in $\Omega$,}
\end{equation}
then
\begin{equation}\label{mphi}
\varphi := m^{\frac{1}{r}} \quad \text{in $\Omega$}
\end{equation}
is a solution of \eqref{rlapl}.
\item[\textbf{b)}] Let $(\varphi, \lambda)$ be a solution of \eqref{rlapl}, and suppose that there exists $u \in C^1(\Omega)$ such that
\begin{equation}\label{uphi2}
h_0 \varphi |Du|^{r'-2}Du + \nu r D \varphi = 0 \quad \text{in $\Omega$.}
\end{equation}
Then, $(u,m,\lambda)$ is a solution of \eqref{MFG}, where $m := \varphi^r$ in $\Omega$.
\end{itemize}
\end{theorem}

The proposed transformation reveals a connection between some (stationary) MFG systems with non-quadratic Hamiltonians of the form \eqref{H} and $r$-Laplace equations, which have been widely studied in the literature and appear in many other areas of interest. Apart from existence and uniqueness issues, this link might shed some light on MFG problems in general, which can be translated into problems involving the $r$-Laplacian (e.g. qualitative properties of solutions, MFG in unbounded domains, vanishing viscosity limit $\nu \rightarrow 0$).

In the quadratic case $r'=2$, we have mentioned that if solutions of \eqref{MFG}-(K) are unique, then $m = e^{-\frac{h_0 u}{\nu}}  / \left( \int_\Omega e^{-\frac{h_0 u}{\nu}} \right)^{-1}$. In this case condition \eqref{crucialeq} is easily verified, and the assertion of Theorem \ref{genHCp} is that $\varphi = m^{1/2}$ solves \eqref{rlapl} with $r=2$, which is precisely \eqref{semilin}. In this sense our transformation can be seen as a generalization of the standard Hopf-Cole.

Note that the change of variables $m = \varphi^r$ is \textit{local}, in particular it is independent of boundary conditions that might be added to \eqref{MFG}. However, in order to verify \eqref{crucialeq}, \eqref{uphi2} and therefore to apply Theorem \ref{genHCp}, it is necessary to specify additional information on the problem. Space dimension $d = 1$ and Neumann conditions at the boundary, or $u,m, \varphi$ enjoying radial symmetry are two possible scenarios where \eqref{crucialeq}, \eqref{uphi2} hold.

\begin{corollary}\label{corradial} Suppose that $H$ satisfies \eqref{H} and $\Omega = \{x \in \Rset^d : |x| < R \}$ for some $R \in (0, \infty]$. Then, $(u,m,\lambda)$ is a radial solution of \eqref{MFG} if and only if $(\varphi, \lambda)$ is a radial solution of \eqref{rlapl}, where $\varphi = m^{1/r}$.
\end{corollary}

\begin{corollary}\label{coroned}Suppose that $H$ satisfies \eqref{H} and $\Omega = (a,b)$ for some $-\infty < a < b < \infty$. Then, $(u,m,\lambda) \in C^2(\overline{\Omega}) \times W^{1,2}(\Omega) \times \Rset$ is a solution of \eqref{MFG} satisfying the Neumann boundary conditions\footnote{For equations \eqref{MFG}-(K) and \eqref{rlapl}, Neumann boundary conditions are intended in the weak sense, namely the space of test functions $\xi$ is set to be $C^\infty(\overline{\Omega})$.} $u'(a) = u'(b) = m'(a) = m'(b) = 0$ if and only if $(\varphi, \lambda) \in W^{1,r}({\Omega}) \cap L^{\infty}({\Omega})$ is a solution of \eqref{rlapl} satisfying the Neumann boundary conditions $\varphi'(a) = \varphi'(b) = 0$, where $\varphi = m^{1/r}$.
\end{corollary}


{\bf Acknowledgements.} The author wishes to express his gratitude to the anonymous referee for his very careful reading of the manuscript and his valuable advices.

\section{Proofs}\label{proofs}

\begin{remark}\label{remregularity} {\normalfont We point out that $m \in W^{1,q}_\loc(\Omega)$ for all $q \ge 1$ by standard regularity results on weak solutions of Kolmogorov equations. Moreover, the Harnack inequality guarantees that $m > 0$ on $\Omega$. However, even if $u \in C^2(\Omega)$, we do not expect in general the same regularity for $m$, since $D H(Du)$ might lack of the desired smoothness if $1 < r' < 2$. If $r' \ge 2$, it is possible to conclude that $m$ is twice differentiable on $\Omega$ and it solves \eqref{MFG}-(K) in the classical sense. \\
On the other hand, it is known that a solution $\varphi$ of \eqref{rlapl} enjoys local $C^{1,\alpha}$ regularity (see, for example, \cite{dibenedetto}, \cite{lieberman}). }
\end{remark}

{\bf Proof of Theorem \ref{genHCp}, a).} Note that $m \in W^{1,q}_\loc(\Omega)$ for all $q \ge 1$ (see Remark \ref{remregularity}), and therefore $\varphi^{r-1}$ has the same regularity. Moreover, $m > 0$ on $\Omega$, hence $\varphi > 0$ as well. Equalities \eqref{crucialeq} and \eqref{mphi} imply that
\begin{equation}\label{uphi}
\nu r \frac{D \varphi}{\varphi} = - h_0 |Du|^{r'-2}Du \quad \text{a.e. in $\Omega$.}
\end{equation}
We multiply the Hamilton-Jacobi-Bellman equation \eqref{MFG}-(HJB) by $\xi \varphi^{r-1}$, where $\xi \in C_0^\infty(\Omega)$ is a generic test function. Integrating by parts,
\begin{equation}\label{hjw}
 \nu \int_{\Omega} D u \cdot D(\xi \varphi^{r-1}) + \int_{\Omega} \frac{h_0}{r'} {|Du(x)|^{r'}} \xi \varphi^{r-1} + \int_{\Omega} (\lambda - f )\xi \varphi^{r-1}= 0.
\end{equation}
We note that by \eqref{uphi}
\[
\nu \int_{\Omega} \xi D u \cdot D(\varphi^{r-1}) = \nu(r-1) \int_{\Omega} \xi D u \cdot D\varphi \, \varphi^{r-2} =
- h_0 \frac{r-1}{r} \int_{\Omega} \xi |D u|^{r'} \varphi^{r-1},
\]
and $(r')^{-1} = (r-1)r^{-1}$, so \eqref{hjw} becomes
\begin{equation}\label{eq1}
\nu \int_{\Omega} D u \cdot D\xi \, \varphi^{r-1} + \int_{\Omega} (\lambda - f )\xi \varphi^{r-1}= 0.
\end{equation}

Again using \eqref{uphi}, we obtain
\begin{multline}\label{weakrlapl}
\int_{\Omega} |D \varphi|^{r-2} D \varphi \cdot D \xi = -\left(\frac{h_0}{\nu r} \right)^{r-1} \int_{\Omega} | D u|^{(r'-1)(r-2) + r' -2} \varphi^{r-1} Du \cdot D \xi \\
= \frac{1}{\nu}\left(\frac{h_0}{\nu r} \right)^{r-1} \int_{\Omega} (\lambda - f )\xi \varphi^{r-1},
\end{multline}
since $(r'-1)(r-2) + r' -2 = 0$, and by virtue of \eqref{eq1}. Equality \eqref{weakrlapl}, which holds for all test functions $\xi$, is precisely the weak formulation of \eqref{rlapl}, hence we are done. We finally observe that $\varphi$ enjoys local $C^{1,\alpha}$ regularity (see Remark \ref{remregularity}), which is inherited by $m$.

{\bf b). } It is easily verified, in view of \eqref{uphi2}, that $m$ solves \eqref{MFG}-(K). Since $\varphi = m^{1/r}$ is positive (see Remark \ref{remregularity}),
\[
\nu r \frac{D \varphi}{\varphi} = - h_0 |Du|^{r'-2}Du \quad \text{in $\Omega$.}
\]
By carrying out backwards computations \eqref{hjw}-\eqref{weakrlapl} of part a), it follows that $u,\lambda$ is a weak solution of \eqref{MFG}-(HJB). Standard regularity results for the Poisson equation guarantee that $u \in C^2(\Omega)$ and \eqref{MFG}-(HJB) is satisfied in the classical sense.
\QED

{\bf Proof of Corollary \ref{corradial}.} In the following, $\rho := |x|$ and $\er = \er(x) = x/\rho$ will denote the standard unit vector in the radial direction.

Let $(u,m,\lambda)$ be a radial solution of \eqref{MFG}. We write $u(x)=u(\rho)$, $m(x)=m(\rho)$, $Du(x)=u'(\rho) \er$, $Dm(x)=m'(\rho) \er$. The Kolmogorov equation \eqref{MFG}-(K) reads
\begin{equation}\label{kolmoexpl}
 \int_\Omega (\nu Dm +  h_0 |Du|^{r'-2} Du \, m) \cdot D \xi = 0
\end{equation}
for all test functions $\xi \in C^{\infty}_0(\Omega)$, and \[\nu Dm(x) +  h_0 |Du(x)|^{r'-2} Du(x) \, m(x) = F(\rho) \er\] for all $x \in \Omega$, where $F(\rho) \rho^{(d-1)/q} \in L^q((0, R'))$, for all $q \ge 1, R' < R$, so \eqref{kolmoexpl} becomes
\[
 \int_0^R F(\rho) \xi'(\rho) \rho^{d-1} d\rho = 0
\]
for all radial test functions $\xi$. It is then possible to conclude that $F(\rho) \rho^{d-1} = 0$ a.e. in $(0, R)$, so \eqref{crucialeq} holds and $(\varphi, \lambda)$ is a solution of \eqref{rlapl} because of Theorem \ref{genHCp}, a). 

If $(\varphi, \lambda)$ is a radial solution of \eqref{rlapl}, say $\varphi(x)=\varphi(\rho)$, $D\varphi(x)=\varphi'(\rho) \er$, then $\varphi'(0)=0$. Setting
\[
b(\rho) := \frac{-\nu r \varphi'(\rho)}{h_0 \varphi(\rho)} \quad \forall \rho \in [0, R)
\]
we have $b \in C^{0,\alpha}([0, R'])$ for all $0 < R' < R$ and $b(0)=0$, so
\[
u(\rho) := \int_0^\rho |b(\sigma)|^{\frac{2-r'}{r'-1}} b(\sigma) d \sigma
\]
defines a radial function $u(x) = u(\rho)$ that belongs to $C^1(\Omega)$. An easy computation shows that $u$ satisfies \eqref{uphi2}, hence $(u,m,\lambda)$ is a solution of \eqref{MFG} as a consequence of Theorem \ref{genHCp}, b).
\QED

{\bf Proof of Corollary \ref{coroned}.} We proceed as in the proof of Corollary \ref{corradial}. Let $(u,m,\lambda)$ be a solution of \eqref{MFG} and set \[\nu m'(x) +  h_0 |u'(x)|^{r'-2} u'(x) \, m(x) =: F(x) \] in $[a,b]$, and $F \in L^q((a,b))$. Therefore,
\[
\int_a^b F(x) \xi'(x) dx = 0
\]
for all test functions $\xi \in C^{\infty}([a,b])$, so $F(x) = 0$ for a. e. $x \in [a,b]$. Hence \eqref{crucialeq} holds and one implication stated by the corollary follows by Theorem \ref{genHCp}, a). Vice-versa, \eqref{uphi2} holds if we choose $u(x) := \int_a^x |b(y)|^{(2-r')/(r'-1)} b(y) d y$, where $b(y) = (-\nu r \varphi'(y))/(h_0 \varphi(y))$ for all $y \in [a,b]$, and Theorem \ref{genHCp}, b) applies.

\QED


\end{document}